\documentclass[12pt]{amsart}
  \usepackage{latexsym}
  \usepackage[all]{xy}
  \usepackage{amsfonts}
  \usepackage{amsthm}
  \usepackage{amsmath}
  \usepackage{amssymb}
  \usepackage[mathscr]{eucal}

 \newcommand{\rmod}[1]{\mathcal{M}_{#1}}
\newcommand{\lmod}[1]{{}_{#1}\mathcal{M}}

\newcommand{\bimod}[2]{{}_{#1}\mathcal{M}_{#2}}

\newcommand{\coring}[1]{\mathfrak{#1}}
\newcommand{\tensor}[1]{\otimes_{#1}}
\newcommand{\tensork}[2]{#1\otimes_k #2}
\newcommand{\tensfun}[1]{\underset{{#1}}{\otimes}}

\newcommand{\bara}[1]{\overline{#1}}
\newcommand{\cotensor}[1]{\square_{#1}}
\newcommand{\cC}{\coring{C}}
\newcommand{\calC}{\mathcal{C}}
\newcommand{\cD}{\coring{D}}
\newcommand{\cE}{\coring{E}}

\newcommand{\CatExt}[2]{\mathbf{CrgExt}^{#1}_{#2}}
\newcommand{\Ecirc}{ \,\, \circ \,\,}

  \newcounter{zlist}

  \newcounter{blist}

  \newcounter{rlist}

  \newtheorem{proposition}{Proposition}[section]
  \newtheorem{lemma}[proposition]{Lemma}

  \theoremstyle{definition}

  \theoremstyle{remark}
  \newtheorem{remark}[proposition]{Remark}

  \newcounter{c}
  \renewcommand{\[}{\setcounter{c}{1}$$}
  \newcommand{\etyk}[1]{\vspace{-7.4mm}$$\begin{equation}\Label{#1}
  \addtocounter{c}{1}}
  \renewcommand{\]}{\ifnum \value{c}=1 $$\else \end{equation}\fi}
\begin{document}

\baselineskip 19pt

  \title{Monoidal categories of corings}
  \author{ Laiachi El Kaoutit}
  \address{Departamento de {\'A}lgebra,  Universidad de
Granada, E18071 Granada, Spain}
\email{kaoutit@ugr.es}
\date{October, 2004}

\begin{abstract}
We introduce a monoidal category of corings using two different
notions of corings morphisms. The first one is the (right) coring
extensions recently introduced by T. Brzezi\'nski in
\cite{Brzezinski:2004}, and the anther is the usual notion of
morphisms defined in \cite{Gomez:2002} by J. G\'omez-Torrecillas.
\end{abstract}

\maketitle

\section*{Introduction}

The word coring appeared for the first time in the literature in
the famous paper of M. Sweedler \cite{Sweedler:1975}, where he
showed that this notion can be used to give a simplest proof of
the first Galois-correspondence theorem for division rings. It
turns out that corings and their comodules unify many kind of
relative modules, like as graded modules, Doi-Hopf modules, and
more generally entwined modules. This was shown in anther famous
paper due to T. Brzezi\'nski \cite{Brzezinski:2002}.

Corings are, in some sense, a generalization of coalgebras to the
case of non-commutative scalars base rings. They have bimodule
structure rather than module one. Thus a tensor product in the
category of bimodules hampers any attempt to give a compatibility
with multiplication and comultiplication. This is a well known
problem to define bialgebras using bimodules. The first approach
to generalize the bialgebras to the case of bimodules was given by
M. Sweedler in \cite{Sweedler:1975a}. If we see that those
bialgebras should be defined as monad or comonad in an appropriate
monoidal category, then we should look first to the possible
monoidal categories. In this note we prove that there is a more
than one monoidal category whose objects are all corings.

We work over a unital commutative ring $k$. All our algebras $A$,
$A'$, $B$, $B'$..., are an unital associative $k$--algebras. For
any algebra we denote by $\rmod{A}$ its category of all unital
right $A$--modules; we use the notation $\lmod{A}$ to denote the
category of unital left $A$--modules. Bimodules are assumed to be
central $k$--bimodules, and their category is denoted by
$\bimod{A}{B}$. If ${}_AM_B$ and ${}_BN_C$ are, respectively, an
$(A,B)$--bimodule and $(B,C)$--bimodule, then their tensor product
$M\tensor{B}N$ will be considered, in a canonical way, as an
$(A,C)$--bimodule.

An $A$--coring is a three-tuple
$(\cC,\Delta_{\cC},\varepsilon_{\cC})$ consisting of an
$A$--bimodule $\cC$ and the two $A$--bilinear maps
$$\xymatrix@C=50pt{\cC \ar@{->}^-{\Delta_{\cC}}[r] & \cC\tensor{A}\cC},\quad
\xymatrix@C=30pt{ \cC \ar@{->}^-{\varepsilon_{\cC}}[r] & A}$$ such
that $(\Delta_{\coring{C}}\tensor{A}\coring{C}) \circ
\Delta_{\coring{C}} = (\coring{C}\tensor{A}\Delta_{\coring{C}})
\circ \Delta_{\coring{C}}$ and
$(\varepsilon_{\coring{C}}\tensor{A}\coring{C}) \circ
\Delta_{\coring{C}}=(\coring{C}\tensor{A}\varepsilon_{\coring{C}})
\circ \Delta_{\coring{C}}= \coring{C}$. A right $\cC$--comodule is
a pair $(M,\rho^{M})$ consisting of right $A$--module and a right
$A$--linear map $\rho^{M}: M \rightarrow M\tensor{A}\cC$, called
right $\cC$--coaction, such that $(M\tensor{A}\Delta_{\coring{C}})
\circ \rho^M = (\rho^M\tensor{A}\coring{C}) \circ \rho^M$ and
$(M\tensor{A}\varepsilon_{\coring{C}}) \circ \rho^M=M$. Left
$\cC$--comodules are symmetrically defined. For instance,
$(\cC,\Delta_{\cC})$ is a left and right $\cC$--comodule. We use
Sweedler's notation for comultiplications, that is
$\Delta_{\cC}(c) = c_{(1)}\tensor{A}c_{(2)}$, for every $c \in
\cC$ (the finite sums are understood). We also use lower indices
Sweedler's notation for coactions: $\rho^M(m) =
m_{(0)}\tensor{A}m_{(1)}$, for all $m\in M $ . But, if two
different coactions are managed, its convenient to use upper
indices: $\rho^{N}(n) = n^{[0]}\tensor{B}n^{[1]}$, for a right
$\cD$--comodule $N$. A source of basic notions on corings,
categories of comodules, bicomodules, and cotensor product is
\cite{Brzezinski/Wisbauer:2003}.

\section{Tensor product of corings}
In this section we recall the tensor product of two corings over
different scalars base rings.

Let $\calC$ and $\calC'$ are, respectively, an $A$--bimodule and
$A'$--bimodule. We consider the tensor product
$\calC\tensor{k}\calC'$ as an $\tensork{A}{A'}$--bimodule using
the canonical bi-action
\begin{eqnarray}\label{biactions}
(\tensork{a}{a'})(\tensork{c}{c'})(\tensork{b}{b'}) &=&
\tensork{(acb)}{(a'c'b')},
\end{eqnarray} for all $(a,b) \in A\times A$, $(a',b')\in A'\times
A'$, and $(c,c')\in \calC\times \calC' $. The following well known
lemma will be used frequently; for completeness we include the
proof.

\begin{lemma}\label{cuatro}
For every pair of modules $(M_A,N_{A'}) \in \rmod{A} \times
\rmod{A'}$, there exists a right $\tensork{A}{A'}$--linear map
\[
\xymatrix@R=0pt{ \eta_{(M_A,N_{A'})}: &
M\tensfun{A}\calC\tensfun{k}N\tensfun{A'}\calC' \ar@{->}[r] &
(\tensork{M}{N}) \tensfun{\tensork{A}{A'}}
(\tensork{\calC}{\calC'})
\\ & m\tensfun{A} c\tensfun{k}n\tensfun{A'}c' \ar@{|->}[r] &
(\tensork{m}{n}) \tensfun{\tensork{A}{A'}} (\tensork{c}{c'}), }
\]
which becomes an $(A,\tensork{A}{A'})$--bilinear map if $M \in \bimod{A}{A}$. Furthermore,
\[
\eta_{(-,-)} : - \tensor{A} \tensork{\calC}{-} \tensor{A'} \calC'
\longrightarrow
(\tensork{-}{-})\tensor{\tensork{A}{A'}}(\tensork{\calC}{\calC'})
\]
is a natural isomorphism.
\end{lemma}
\begin{proof}
It is clear that $A\tensor{A}\calC\tensor{k}A'\tensor{A'} \calC'
\cong
(\tensork{A}{A'})\tensor{\tensork{A}{A'}}(\tensork{\calC}{\calC'})$,
via the map sending $a\tensor{A}c\tensor{k}a'\tensor{A'}c' \mapsto
(\tensork{a}{a'}).(\tensork{c}{c'})=\tensork{ac}{a'c'}$. If
$(f,g): (A_A,A'_{A'}) \rightarrow (A_A,A'_{A'})$ is any arrow in
the product category $\rmod{A} \times \rmod{A'}$, then its also
clear that
\[
\xymatrix@R=30pt@C=70pt{
A\tensfun{A}\calC\tensfun{k}A'\tensfun{A'} \calC'
\ar@{->}^-{f\tensfun{A}\calC\tensfun{k}g\tensfun{A'} \calC'}[r]
\ar@{->}^-{\cong}[d] &
A\tensfun{A}\calC\tensfun{k}A'\tensfun{A'} \calC'  \ar@{->}^-{\cong}[d] \\
(\tensork{A}{A'})\tensfun{\tensork{A}{A'}}(\tensork{\calC}{\calC'})
\ar@{->}^-{(\tensork{f}{g})\tensfun{\tensork{A}{A'}}(\tensork{\calC}{\calC'})}[r]
& (\tensork{A}{A'}) \tensfun{\tensork{A}{A'}}
(\tensork{\calC}{\calC'}) }
\]
is a commutative diagram. Since $(A_A,A'_{A'})$ is a projective
generator in $\rmod{A} \times \rmod{A'}$ and the tensor product
commute with direct limits, Mitchell's Theorem \cite[Theorem 5.4,
p. 109]{Mitchell:1965}, implies that there exists a unique natural
isomorphism
\[
\eta_{(-,-)} :  - \tensor{A} \tensork{\calC}{-} \tensor{A'} \calC'
\longrightarrow
(\tensork{-}{-})\tensor{\tensork{A}{A'}}(\tensork{\calC}{\calC'}).
\]
Let $m \in M$, $n \in N$, $c\in \calC$ and $c' \in \calC'$, the
form of the image
$\eta_{(M,N)}(m\tensor{A}c\tensor{k}n\tensor{A'}c')$, is computed
by using the morphisms $(A_A,A'_{A'}) \rightarrow (mA,nA')$ and
the naturality of $\eta$, where $mA$ and $nA$ denote the cyclic
submodules.
\end{proof}

Let $(\cC,\Delta_{\cC},\varepsilon_{\cC})$ and $(\cC',\Delta_{\cC'},\varepsilon_{\cC'})$ are, respectively,
an $A$--coring and $A'$--coring, and consider $\tensork{\cC}{\cC'}$ canonically as an $\tensork{A}{A'}$--bimodule.
\begin{proposition}[\cite{Gomez/Louly:2003}]\label{tens-corg}
The tensor product $\cC\otimes_k \cC'$ is an  $A\otimes_k
A'$--coring with comultiplication  given by the composition map
$$\xymatrix@C=50pt{ \cC\otimes_k \cC'
\ar@{->}^-{\Delta_{\cC}\otimes_k\Delta_{\cC'}}[r] &
(\cC\otimes_A\cC)\otimes_k (\cC'\otimes_{A'}\cC')
\ar@{->}^-{\eta_{\cC,\cC'}}[r] &
(\cC\otimes_k\cC')\otimes_{A\otimes_k A'} (\cC\otimes_{k}\cC'),
}$$ and counit by $$\xymatrix@C=50pt{ \tensork{\cC}{\cC'}
\ar@{->}^-{\tensork{\varepsilon_{\cC}}{\varepsilon_{\cC'}}}[r] &
\tensork{A}{A'}. }$$
\end{proposition}

\section{Tensor product of right corings extensions}

Let $A$ and $B$ be $k$--algebras, $\cC$ an $A$--coring and $\cD$ a $B$--coring.
Recall from \cite[Definition 2.1]{Brzezinski:2004},
that $\cD$ is called \emph{right extension} of $\cC$ provided $\cC$ is a $(\cC,\cD)$--bicomodule
with the regular left coaction $\Delta_{\cC}$. This means that $\cC$ should be an
$(A,B)$--bimodule and $\Delta_{\cC}$ a right $B$--linear map.

\begin{proposition}\label{tensor}
Let $\cD$ and $\cD'$ are, respectively, a $B$--coring and $B'$--coring. Assume that $\cD$
(resp. $\cD'$) is a right extension of $\cC$ (resp. of $\cC'$). Then $\cD\otimes_k\cD'$ is a
right extension of $\cC\otimes_k\cC'$.

\end{proposition}
\begin{proof}
Denote by $\rho^{\cC}: \cC \rightarrow \cC\otimes_B\cD$ ($c\mapsto
c_{(0)}\otimes_Bc_{(1)}$) and $\rho^{\cC'}: \cC' \rightarrow
\cC'\otimes_{B'}\cD'$ ($c'\mapsto c'_{(0)}\otimes_{B'}c'_{(1)}$)
(sums are understood) the attached right coaction to the stated
extensions. Define $\rho^{\cC\otimes_k\cC'}$ as the composition
map $$\xymatrix@R=0pt@C=40pt{\cC\otimes_k\cC'
\ar@{->}^-{\rho^{\cC}\otimes_k \rho^{\cC'}}[r] &
(\cC\otimes_B\cD)\otimes_k(\cC'\otimes_{B'}\cD')
\ar@{->}^-{\eta_{\cC,\cC'}}[r]&
(\cC\otimes_k\cC')\otimes_{B\otimes_kB'}(\cD\otimes_k\cD') \\
c\otimes_kc' \ar@{|->}[rr] & &  (c_{(0)}\otimes_k
c'_{(0)})\otimes_{B\otimes_kB'}(c_{(1)}\otimes_kc'_{(1)})}$$ It is
easily checked, using the $B$--linearity of $\Delta_{\cC}$ and the
$B'$--linearity of $\Delta_{\cC'}$, that this composition is
$(A\otimes_k A')-(B\otimes_k B')$--bilinear map. Since
$\rho^{\cC}$ and $\rho^{\cC'}$ are, respectively, right
$\cD$--coaction and right $\cD'$--coaction and $\eta_{-,-}$ is a
natural transformation, $\rho^{\cC\otimes_k\cC'}$ is also a right
$\cD\otimes_k\cD'$--coaction. Furthermore, $\rho^{\cC}$,
$\rho^{\cC'}$, and the comultiplications of $\cC$ and $\cC'$,
enjoy the following four commutative diagrams

$$\xymatrix@R=40pt@C=80pt{\cC\tensfun{k}\cC' \ar@{->}|-{\Delta\tensfun{k}\Delta'}[d]
\ar@{->}^-{\rho^{\cC}\tensfun{k}\rho^{\cC'}}[r] \ar@{}^-{(1)}[rd] & (\cC\tensfun{B}\cD)
\tensfun{k}(\cC'\tensfun{B'}\cD')  \ar@{->}|-{(\Delta\otimes_B\cD)\otimes_k(\Delta'\otimes_{B'}\cD')}[d]
\\ (\cC\tensfun{A}\cC)\tensfun{k}(\cC'\tensfun{A'}\cC') \ar@{->}^-{(\cC\tensfun{A}\rho^{\cC})
\tensfun{k}(\cC'\tensfun{A'}\rho^{\cC'})}[r]  &   (\cC\tensfun{A}\cC\tensfun{B}\cD)\tensfun{k}
(\cC'\tensfun{A'}\cC'\tensfun{B'}\cD')
}$$

$$\xymatrix@R=40pt@C=45pt{ (\cC\tensfun{B}\cD)\tensfun{k}(\cC'\tensfun{B'}\cD')
\ar@{->}^-{\cong}[r] \ar@{->}|-{(\Delta\otimes_B\cD)\otimes_k(\Delta'\otimes_{B'}\cD')}[d]
\ar@{}^-{(2)}[rd] &
(\cC\tensfun{k}\cC')\tensfun{B\otimes_kB'}(\cD\tensfun{k}\cD') \ar@{->}|-{(\Delta\tensfun{k}\Delta')
\tensfun{B\otimes_kB'}(\cD\tensfun{k}\cD')}[d]
\\ (\cC\tensfun{A}\cC\tensfun{B}\cD)\tensfun{k}(\cC'\tensfun{A'}\cC'\tensfun{B'}\cD') \ar@{->}^-{\cong}[r]
& \left( (\cC\tensfun{A}\cC)\tensfun{k}(\cC'\tensfun{A'}\cC')\right)\tensfun{B\otimes_kB'}(\cD\tensfun{k}\cD')   }$$

$$\xymatrix@R=40pt@C=50pt{ (\cC\tensfun{A}\cC)\tensfun{k}(\cC'\tensfun{A'}\cC')
\ar@{->}^-{(\cC\tensfun{A}\rho^{\cC})\tensfun{k}(\cC'\tensfun{A'}\rho^{\cC'})}[r]
\ar@{->}|-{\cong}[d] \ar@{}^-{(3)}[rd] &
(\cC\tensfun{A}\cC\tensfun{B}\cD)\tensfun{k}(\cC'\tensfun{A'}\cC'\tensfun{B'}\cD')
\ar@{->}|-{\cong}[d]
\\ (\cC\tensfun{k}\cC')\tensfun{A\otimes_kA'}(\cC\tensfun{k}\cC')
\ar@{->}^-{(\cC\tensfun{k}\cC')\tensfun{A\otimes_kA'}(\rho^{\cC}\tensfun{k}\rho^{\cC'})}[r]
&  (\cC\tensfun{k}\cC') \tensfun{A\otimes_kA'}\left( (\cC\otimes_{B}\cD)
 \tensfun{k} (\cC'\otimes_{B'}\cD') \right)  }$$

$$\xymatrix@R=40pt@C=20pt{   (\cC\tensfun{A}\cC\tensfun{B}\cD)\tensfun{k}(\cC'\tensfun{A'}\cC'\tensfun{B'}\cD')
 \ar@{->}^-{\cong}[r] \ar@{->}|-{\cong}[d] \ar@{}^-{(4)}[rd] & \left( (\cC\tensfun{A}\cC)\tensfun{k}
 (\cC'\tensfun{A'}\cC')\right)\tensfun{B\otimes_kB'}(\cD\tensfun{k}\cD') \ar@{->}|-{\cong}[d]
\\  (\cC\tensfun{k}\cC') \tensfun{A\otimes_kA'}\left( (\cC\otimes_{B}\cD) \tensfun{k} (\cC'\otimes_{B'}\cD') \right)
 \ar@{->}^-{\cong}[r] & (\cC\tensfun{k}\cC') \tensfun{A\otimes_kA'} (\cC\tensfun{k}\cC') \tensfun{B\otimes_k B'}
 (\cD\tensfun{k}\cD') }$$ where the isomorphisms maps are defined by the natural isomorphism of Lemma \ref{cuatro}.
If we put those diagrams in the following form
$$\begin{tabular}{|c|c|}
\hline  (1) & (2) \\
\hline  (3) & (4) \\  \hline \end{tabular}$$ we then get another
commutative diagram which shows that $\rho^{\cC\otimes_k\cC'}$ is
left $\cC\otimes_{k}\cC'$--colinear with respect to regular left
coaction $\Delta_{\tensork{\cC}{\cC'}}$, and this finishes the
proof.
\end{proof}

\section{A monoidal category}

Let us recall from \cite{Brzezinski:2004} the category of corings
$\CatExt{r}{k}$. The object of this category are corings
understood as pair $(\cC:A)$ (that is $\cC$ is an $A$--coring),
and morphisms $(\cC:A) \rightarrow (\cD:B)$ are pairs
$(\rho_{\cC},\rho^{\cC})$ where $\rho_{\cC}: \cC\otimes B
\rightarrow \cC$ is left $\cC$--colinear right $B$--action, and
$\rho^{\cC}: \cC \rightarrow \cC\tensor{B}\cD$ is a left
$\cC$--colinear right $\cD$--coaction (that is $\cD$ is right
extension of $\cC$). The identity arrow of an object $(\cC:A)$ is
given by the pair $id_{(\cC:A)}=(\rho_{\cC},\rho^{\cC})$ where
$\rho_{\cC}=\iota_{\cC}: \cC \tensor{k} A \rightarrow \cC$ is the
initial right $A$--action and $\rho^{\cC}=\Delta_{\cC}$ is the
comultiplication of $\cC$. The composition law is given as
follows. If $(\rho_{\cE},\rho^{\cE}): (\cE:C)\rightarrow (\cC:A)$
and $(\rho_{\cC},\rho^{\cC}): (\cC:A) \rightarrow (\cD:B)$, then
$$(\rho_{\cC},\rho^{\cC}) \circ (\rho_{\cE},\rho^{\cE}) \, = \,
(\rho_{\cC}\bullet \rho_{\cE}, \rho^{\cC} \bullet \rho^{\cE}), $$
where $$\xymatrix@C=40pt{ \rho_{\cC}\bullet \rho_{\cE}:
\cE\otimes_k B \ar@{->}^{\rho^{\cE}\otimes B}[r] & \cE
\tensor{A}\cC \tensor{k} B \ar@{->}^-{\cE\tensor{A}\rho_{\cC}}[r]
& \cE\tensor{A}\cC \ar@{->}^-{\cE\tensor{A}\varepsilon}[r] &
\cE\tensor{A}A \cong \cE }$$ and $$\xymatrix@C=60pt{
\rho^{\cC}\bullet \rho^{\cE}: \cE \ar@{->}^-{\cong}[r] &
\cE\cotensor{\cC}\cC \ar@{->}^-{\cE\cotensor{\cC}\rho^{\cC}}[r] &
\cE \cotensor{\cC}(\cC \tensor{B}\cD) \ar@{->}^-{\cong}[r] &
\cE\tensor{B} \cD.  }$$ Explicitly, the bullet compositions are
given as follows: for $e \in\cE$ and $b \in B$ \begin{eqnarray*}
 \rho_{\cC}\bullet \rho_{\cE} (e\tensor{k}b) &=&  e_{(0)}\varepsilon_{\cC}(e_{(1)} b), \\
\rho^{\cC}\bullet \rho^{\cE} (e) &=& e_{(0)} \varepsilon_{\cC}(e_{(1)}{}^{[0]} ) \tensor{B} e_{(1)}{}^{[1]}
\end{eqnarray*}
where $\rho^{\cE}(e) = e_{(0)} \tensor{A}e_{(1)}$, and
$\rho^{\cC}(e_{(1)}) = e_{(1)}{}^{[0]} \tensor{B}
e_{(1)}{}^{[1]}$, for all $e_{(1)}$.

The tensor product of two morphisms
$(\rho_{\cC},\rho^{\cC}):(\cC:A) \rightarrow (\cD:B)$ and
$(\rho_{\cC'},\rho^{\cC'}): (\cC':A') \rightarrow (\cD':B')$ in
$\CatExt{r}{k}$, is defined as in the proof of Proposition
\ref{tensor}; that is by the following morphism:
\begin{equation}\label{tens-map}
(\rho_{\tensork{\cC}{\cC'}},\rho^{\tensork{\cC}{\cC'}}):
(\tensork{\cC}{\cC'}: \tensork{A}{A'}) \longrightarrow
(\tensork{\cD}{\cD'}:\tensork{B}{B'}) \end{equation}
$$\xymatrix@R=0pt{\rho_{\tensork{\cC}{\cC'}}: (\tensork{\cC}{\cC'}
\tensor{k} (\tensork{B}{B'}) \ar@{->}[r] & \tensork{\cC}{\cC'}
\\  (\tensork{c}{c'}) \tensor{k} (\tensork{b}{b'}) \ar@{|->}[r]
& (cb)\tensor{k} (c'b')  }$$ where $\rho_{\cC}(c\tensor{k}b)=cb$,
$\rho_{\cC'}(c'\tensor{k}b')=c'b'$, and
$$\xymatrix@R=0pt{\rho^{\tensork{\cC}{\cC'}}: \cC\tensor{k}\cC'
\ar@{->}[r] & (\cC\tensor{k}\cC') \tensor{\tensork{B}{B'}}
(\tensork{\cD}{\cD'} \\ c\tensor{k}c' \ar@{|->}[r] &
(c_{(0)}\tensor{k}c'_{(0)}) \tensor{\tensork{B}{B'}}
(c_{(1)}\tensor{k}c'_{(1)})  }$$ where $\rho^{\cC}(c) =
c_{(0)}\tensor{B}c_{(1)}$, $\rho^{\cC'}(c') =
c'_{(0)}\tensor{B'}c'_{(1)}$.

\begin{proposition}\label{bifunt}
Let $k$ be an unital commutative ring. Consider the category
$\CatExt{r}{k}$ of corings with morphisms right coring extensions,
and denote by $\Bbbk:=(k:k)$ the trivial $k$--coring $k$. There
exists a covariant bi-functor
$$\xymatrix@R=0pt{-\tensor{\Bbbk}-: \CatExt{r}{k} \times
\CatExt{r}{k} \ar@{->}[r] & \CatExt{r}{k} \\ \left(
(\cC:A),(\cC':A') \right) \ar@{->}[r] &  (\tensork{\cC}{\cC'}:
\tensork{A}{A'}) \\ \left(
(\rho_{\cC},\rho^{\cC}),(\rho_{\cC'},\rho^{\cC'})\right)
\ar@{->}[r] & (\rho_{\tensork{\cC}{\cC'}},
\rho^{\tensork{\cC}{\cC'}}) }$$ where
$(\rho_{\tensork{\cC}{\cC'}}, \rho^{\tensork{\cC}{\cC'}})$ is the
morphism defined in equation \eqref{tens-map}. Moreover, $$\Bbbk
\tensor{\Bbbk} (\cC:A) \cong (\cC:A) \text{ and }
(\cC:A)\tensor{\Bbbk} \Bbbk \cong (\cC:A)$$ natural isomorphisms
in $\CatExt{r}{k}$. In particular, $\CatExt{r}{k}$ is a monoidal
category with unit $\Bbbk$.
\end{proposition}
\begin{proof}
After all, the Propositions \ref{tens-corg} and \ref{tensor},
imply that the stated functor is well defined. Now, by definition
of the comultiplication of the tensor product of two corings, the
identity arrow of any pair of corings is mapped by
$-\tensor{\Bbbk}-$ to the identity arrow of their tensor product;
that is in the above notation, we have
$$id_{(\cC:A)} \tensor{\Bbbk}id_{(\cC':A')}
\,=\,(\iota_{\cC},\Delta_{\cC})\tensor{\Bbbk}
(\iota_{\cC'},\Delta_{\cC'}) \,=\,
(\iota_{\cC\tensor{k}\cC'},\Delta_{\tensork{\cC}{\cC'}}) \,=\,
id_{ (\cC\tensor{k}\cC':A\tensor{k}A')}.$$ Consider the following
four morphisms in $\CatExt{r}{k}$ $$\xymatrix@C=40pt{ (\cE:C)
\ar@{->}^-{(\rho_{\cE},\rho^{\cE})}[r] & (\cC:A)
\ar@{->}^-{(\rho_{\cC},\rho^{\cC})}[r] & (\cD:B) \\  (\cE':C')
\ar@{->}^-{(\rho_{\cE'},\rho^{\cE'})}[r] & (\cC':A')
\ar@{->}^-{(\rho_{\cC'},\rho^{\cC'})}[r] & (\cD':B'),  }$$ and put
$(\rho_{\cC}\bullet \rho_{\cE}, \rho^{\cC} \bullet \rho^{\cE})
\tensor{\Bbbk} (\rho_{\cC'}\bullet \rho_{\cE'}, \rho^{\cC'}
\bullet \rho^{\cE'})
=(\bara{\rho}_{\tensork{\cE}{\cE'}},\bara{\rho}^{\tensork{\cE}{\cE'}})$.
By definition $\bara{\rho}_{\tensork{\cE}{\cE'}}:
(\tensork{\cE}{\cE'}) \tensor{k} (\tensork{B}{B'}) \rightarrow
\tensork{\cE}{\cE'}$ sends $(\tensork{e}{e'}) \tensor{k}
(\tensork{b}{b'}) \mapsto (eb)\tensor{k}(e'b')$, where $eb=
\rho_{\cC} \bullet \rho_{\cE} (e\tensor{k}b) =
e_{(0)}\varepsilon_{\cC}(e_{(1)}b)$ and  $e'b'= \rho_{\cC'}
\bullet \rho_{\cE'} (e'\tensor{k}b') =
e'_{(0)}\varepsilon_{\cC'}(e'_{(1)}b')$. That is \begin{eqnarray*}
\bara{\rho}_{\tensork{\cE}{\cE'}} \left( (\tensork{e}{e'})
\tensor{k} (\tensork{b}{b'}) \right) &=&
e_{(0)}\varepsilon_{\cC}(e_{(1)}b) \tensor{k}
e'_{(0)}\varepsilon_{\cC'}(e'_{(1)}b') \\ &=&
\rho_{\tensork{\cC}{\cC'}} \bullet \rho_{\tensork{\cE}{\cE'}}
\left( (\tensork{e}{e'}) \tensor{k} (\tensork{b}{b'}) \right)
\end{eqnarray*} for every $e \in \cE$, $e'\in \cE$, $b \in B$,
and $b' \in B'$. Thus $\bara{\rho}_{\tensork{\cE}{\cE'}} =
\rho_{\tensork{\cC}{\cC'}} \bullet \rho_{\tensork{\cE}{\cE'}}$.

From another hand the map $\bara{\rho}^{\tensork{\cE}{\cE'}}$ is defined by the composition
$$\xymatrix@R=0pt@C=50pt{ \cE\tensor{k}\cE' \ar@{->}^-{(\rho^{\cC} \bullet \rho^{\cE})
\otimes(\rho^{\cC'} \bullet \rho^{\cE'})}[r] & (\cE\tensor{B}\cD)
\tensor{k} (\cE'\tensor{B'}\cD') \ar@{->}^-{\eta_{\cE,\cE'}}[r] &
(\tensork{\cE}{\cE'}) \tensor{\tensork{B}{B'}}
(\tensork{\cD}{\cD'}) } $$ sending $$ \xymatrix@C=50pt{
\tensork{e}{e'} \ar@{|->}[r] & \left(
(e_{(0)}\varepsilon_{\cC}(e_{(1)}{}^{[0]})) \tensor{k} (e'_{(0)}
\varepsilon_{\cC'}(e'_{(1)}{}^{[0]}) ) \right)
\tensor{\tensork{B}{B'}} \left( e_{(1)}{}^{[1]} \tensor{k}
e'_{(1)}{}^{[1]} \right),  }$$ that is
$\bara{\rho}^{\tensork{\cE}{\cE'}}(\tensork{e}{e'}) =
\rho^{\tensork{\cC}{\cC'}} \bullet \rho^{\tensork{\cE}{\cE'}}
(\tensork{e}{e'})$, for every $(e,e') \in \cE \times \cE'$. Hence
$\bara{\rho}^{\tensork{\cE}{\cE'}}= \rho^{\tensork{\cC}{\cC'}}
\bullet \rho^{\tensork{\cE}{\cE'}}$. Therefore, \begin{eqnarray*}
(\rho_{\cC}\bullet \rho_{\cE}, \rho^{\cC} \bullet \rho^{\cE})
\tensor{\Bbbk} (\rho_{\cC'}\bullet \rho_{\cE'}, \rho^{\cC'}
\bullet \rho^{\cE'}) \\ \, =\, (\rho_{\tensork{\cC}{\cC'}} \bullet
\rho_{\tensork{\cE}{\cE'}},\rho^{\tensork{\cC}{\cC'}} \bullet
\rho^{\tensork{\cE}{\cE'}}) \\ \,=\,
(\rho_{\cC},\rho^{\cC})\tensor{\Bbbk} (\rho_{\cC'},\rho^{\cC'}) \,
\circ \, (\rho_{\cE},\rho^{\cE})\tensor{\Bbbk}
(\rho_{\cE'},\rho^{\cE'}) \\  \,=\, \left( (\rho_{\cC},\rho^{\cC})
\circ (\rho_{\cE},\rho^{\cE}) \right) \tensor{\Bbbk} \left(
(\rho_{\cC'},\rho^{\cC'}) \circ (\rho_{\cE'},\rho^{\cE'}) \right)
\end{eqnarray*} which shows that $-\tensor{\Bbbk}-$ is a covariant
functor. The last assertion is obvious, and the  particular one is
clear, since the associativity of the bi-functor
$-\tensor{\Bbbk}-$ is up to natural isomorphisms.
\end{proof}

\begin{remark}\label{symmetric}
Of course, we have a similar results on the left hand. That is,
the category $\CatExt{l}{k}$ whose object are corings and
morphisms are left coring extensions is also a monoidal category
with the same unit $\Bbbk=(k:k)$. If we denote by $\CatExt{}{k}$
the category whose object are all corings and morphisms are the
left and right (at the same times) coring extensions. Then the
study of $\CatExt{}{k}$ can be also posed, viewing it as a
subcategory of both $\CatExt{r}{k}$ and $\CatExt{l}{k}$. Examples
of morphisms in this subcategory are morphisms between corings
with the same scalars base ring.
\end{remark}

Finally, we will consider the category of corings
$\coring{Corings}$ where the objects are corings understood also
as pairs $(\cC:A)$ and morphisms are in the sense of
\cite{Gomez:2002}; that is a pair of maps $(\phi,\varphi): (\cC:A)
\rightarrow (\cD:B)$, where $\varphi: A \rightarrow B$ is an
algebra map and $\phi: \cC \rightarrow \cD$ is an $A$--bilinear
map ($\cD$ is an $A$--bimodule by scalars restriction) such that
the following diagrams commute
\[
\xymatrix@C=60pt{\coring{C} \ar@{->}^{\Delta_{\coring{C}}}[r]
\ar@{->}_{\phi}[dd]& \coring{C} \tensor{A} \coring{C}
\ar@{->}^{\phi \tensor{A} \phi }[rd] &  \\ & & \coring{D}
\tensor{A} \coring{D}, \ar@{->}^{\omega_{A,B}}[ld]
\\ \coring{D} \ar@{->}^{\Delta_{\coring{D}}}[r] & \coring{D}
\tensor{B} \coring{D} & } \hspace{2em} \xymatrix@C=45pt{\coring{C}
\ar@{->}^{\varepsilon_{\coring{C}}}[r]\ar@{->}_{\phi}[dd] & A
\ar@{->}^{\varphi}[dd]
\\ & \\ \coring{D} \ar@{->}^{\varepsilon_{\coring{D}}}[r] & B, }
\]
where $\omega_{A,B} :\coring{D} \tensor{A} \coring{D} \rightarrow
\coring{D} \tensor{B} \coring{D}$ is the obvious map associated to
$\varphi$. The identity arrow of an object $(\cC:A)$ is
$id_{(\cC:A)}=(id_{\cC},id_A)$, and the composition law is
componentwise composition.

\begin{lemma}\label{tensor-morph}
Let $(\phi,\varphi): (\cC:A) \rightarrow (\cD:B)$ and
$(\phi',\varphi'): (\cC':A') \rightarrow (\cD':B')$ are a corings
morphisms. Then
$$(\tensork{\phi}{\phi'},\varphi\tensor{k}\varphi'):
(\tensork{\cC}{\cC'}:\tensork{A}{A'}) \rightarrow
(\tensork{\cD}{\cD'},\tensork{B}{B'})$$ is also a corings
morphism.
\end{lemma}
\begin{proof}
Analogue to that of  \cite[Proposici\'on 1.1.20]{Kaoutit:PhD}.
Obviously $\tensork{\varphi}{\varphi'}$ is an algebra map. Since
$\phi$ and $\phi'$ are, respectively, $A$--bilinear and
$A'$--bilinear, it is easily checked that $\phi\tensor{k}\phi'$ is
$\tensork{A}{A'}$--bilinear. The counit property of
$\phi\tensor{k}\phi'$ is given by the following commutative
diagram $$\xymatrix@C=45pt{\coring{C} \tensor{k}\cC'
\ar@{->}^{\varepsilon_{\coring{C}}\tensor{k}\varepsilon_{\cC'}}[r]\ar@{->}_{\phi\tensor{k}\phi'}[dd]
& A\tensor{k}A' \ar@{->}^{\varphi\tensor{k}\varphi'}[dd]
\\ & \\ \coring{D}\tensor{k}\cD' \ar@{->}^{\varepsilon_{\coring{D}}
\tensor{k}\varepsilon_{\coring{D}'}}[r] & B\tensor{k}B'. }$$
For the coassociativity, we denote by $\eta_{-,-}^{A,A'}$ and
$\eta_{-,-}^{B,B'}$ the natural isomorphism of Lemma \ref{cuatro}
to distinguish between the tensor product algebra
$\tensork{A}{A'}$ and $\tensork{B}{B'}$. With this notation, we
have \begin{eqnarray*} \omega_{\tensork{A}{A'},\tensork{B}{B'}}
\Ecirc \left( (\tensork{\phi}{\phi'}) \tensor{\tensork{A}{A'}}
(\tensork{\phi}{\phi'}) \right) \Ecirc
\Delta_{\tensork{\cC}{\cC'}} \\ \,=\,
\omega_{\tensork{A}{A'},\tensork{B}{B'}} \Ecirc \left(
(\tensork{\phi}{\phi'}) \tensor{\tensork{A}{A'}}
(\tensork{\phi}{\phi'}) \right) \Ecirc \eta_{\cC,\cC'}^{A,A'}
\Ecirc (\Delta_{\cC}\tensor{k}\Delta_{\cC'}) \\ \,=\,
\omega_{\tensork{A}{A'},\tensork{B}{B'}} \Ecirc
\eta_{\cD,\cD'}^{A,A'} \Ecirc   \left(
(\phi\tensor{A}\phi)\tensor{k} (\phi'\tensor{A'}\phi') \right)
\Ecirc (\Delta_{\cC}\tensor{k}\Delta_{\cC'}) \\ \,=\,
\eta_{\cD,\cD'}^{B,B'} \Ecirc \left( \omega_{A,B} \tensor{k}
\omega_{A',B'}\right)  \Ecirc   \left(
(\phi\tensor{A}\phi)\tensor{k} (\phi'\tensor{A'}\phi') \right)
\Ecirc (\Delta_{\cC}\tensor{k}\Delta_{\cC'}) \\ \, = \,
\eta_{\cD,\cD'}^{B,B'} \Ecirc
(\Delta_{\cD}\tensor{k}\Delta_{\cD'}) \Ecirc (\phi\tensor{k}\phi')
\,=\, \Delta_{\cD\tensor{k}\cD'} \Ecirc (\phi\tensor{k}\phi'),
\end{eqnarray*} where we have used the naturality of
$\eta_{-,-}^{-,-}$ and the coassociativity of $\phi$ and $\phi'$;
this means the coassociative property of the proposed map.

\end{proof}

\begin{proposition}
Consider the category $\coring{Corings}$ of all corings, and
denote by $\Bbbk=(k:k)$ the trivial $k$--coring. The following
$$\xymatrix@R=0pt{-\tensor{\Bbbk}-: \coring{Corings} \times
\coring{Corings} \ar@{->}[r] & \coring{Corings} \\ \left(
(\cC:A),(\cC':A') \right) \ar@{->}[r] &  (\tensork{\cC}{\cC'}:
\tensork{A}{A'}) \\ \left( (\phi,\varphi),(\phi',\varphi')\right)
\ar@{->}[r] & (\tensork{\phi}{\phi'}, \tensork{\varphi}{\varphi'})
}$$ establishes a covariant bi-functor. Moreover, $$\Bbbk
\tensor{\Bbbk} (\cC:A) \cong (\cC:A) \text{ and }
(\cC:A)\tensor{\Bbbk} \Bbbk \cong (\cC:A)$$ natural isomorphisms
in $\coring{Corings}$. In particular, $\coring{Corings}$ is a
monoidal category with unit $\Bbbk$.

\end{proposition}
\begin{proof}
Consequence of Lemma \ref{tensor-morph}.
\end{proof}

\begin{remark}
The relationship between the morphisms of $\coring{Corings}$ and
those of $\CatExt{r}{k}$ can be given as follows. Let
$(\phi,\varphi): (\cC:A) \rightarrow (\cD:B)$ be a morphism in
$\coring{Corings}$, and consider the $B$--coring
$B\tensor{A}\cC\tensor{A}B$ called the base ring extension of
$\cC$, see \cite[17.2]{Brzezinski/Wisbauer:2003}. Denote by
$\rho_{B\tensor{A}\cC\tensor{A}B}:
B\tensor{A}\cC\tensor{A}B\tensor{k}B \rightarrow
B\tensor{A}\cC\tensor{A}B$ the right scalar multiplication map,
and define $$\xymatrix@R=0pt{ \rho^{B\tensor{A}\cC\tensor{A}B}:
B\tensor{A}\cC\tensor{A}B \ar@{->}[r] & B\tensor{A}\cC\tensor{A}B
\tensor{B}\cD \cong B\tensor{A}\cC\tensor{A}\cD \\
b\tensor{A}c\tensor{A}b' \ar@{|->}[r] & b\tensor{A}c_{(1)}
\tensor{A} \phi(c_{(2)})b' }$$ where
$\Delta_{\cC}(c)=c_{(1)}\tensor{A}c_{(2)}$. It is easily checked
that the pair
$$ (\rho_{B\tensor{A}\cC\tensor{A}B}\, ,\, \rho^{B\tensor{A}\cC\tensor{A}B}):
(B\tensor{A}\cC\tensor{A}B:B) \longrightarrow (\cD:B) $$ is now a
morphism in the category $\CatExt{r}{k}$.

In fact $B\tensor{A}\cC\tensor{A}B$ is a right and left coring
extensions of $\cD$. That is to any morphism $(\phi,\varphi)$ in
$\coring{Corings}$, one can associate a morphism in the category
$\CatExt{}{k}$ described in Remark \ref{symmetric}.
\end{remark}

\section*{Acknowledgements}

The author started working on this note during his visit to the
Mathematics Laboratory (CNRS) of the University of Reims. He would
like to thank Jacques Alev, St\'ephane Launois and other members
of the Laboratory for very warm hospitality. This research has
been partly developed during my postdoctoral stay at the
Universidad de Granada supported by the grant SB2003-0064 from the
Mi\-nis\-terio de Educaci{\'o}n, Cultura y Deporte of Spain.

\end{document}